\documentclass[11pt]{article} 

\usepackage[utf8]{inputenc}
\usepackage{caption,graphicx,amsmath,amsthm,amssymb, authblk,xcolor,enumitem,float,verbatim,hyperref,pgfplots,bibentry,mathtools,soul,isomath,multicol,wrapfig}
\usepackage{mathrsfs}
\usepackage[utf8]{inputenc}
\usepackage{fullpage}
\usepackage{amsfonts, theoremref, authblk, subcaption, multirow, lipsum, booktabs, array, longtable,etoolbox}
\usepackage{bigdelim}
\usepackage{blkarray}
\usepackage{indentfirst}
\usepackage{nccmath}
\graphicspath{ {images/} }
\usepackage{multicol}
\usepackage{graphicx,amsmath,amssymb,amsthm,marvosym,tikz,fontenc}
\usetikzlibrary{backgrounds}
\usetikzlibrary{patterns,shapes.misc, positioning}
\usepackage{pgfplots}
\pgfplotsset{width=10cm,compat=1.9,tick scale binop=\times}
\usepackage{indentfirst}
\usepackage[a4paper,bindingoffset=0.2in,%
left=0.6in,right=0.8in,top=1.2in,bottom=1.2in,%
footskip=.25in]{geometry}
\usepackage{relsize}
\usepackage[english]{babel}
\allowdisplaybreaks
\theoremstyle{plain}
\newtheorem{theorem}{Theorem}[section]
\newtheorem{lemma}[theorem]{Lemma}

\newtheorem{corollary}[theorem]{Corollary}
\newtheorem{definition}[theorem]{Definition}

\newtheorem{example}[theorem]{Example}
\newtheorem{thm}{Theorem}[subsection]
\newtheorem{defn}[thm]{Definition}
\newtheorem{corlly}[thm]{Corollary}
\newtheorem{exmp}[thm]{Example}
\DeclareRobustCommand{\rchi}{{\mathpalette\irchi\relax}} 
\newcommand{\irchi}[2]{\raisebox{\depth}{$#1\chi$}}   

\newlength{\defbaselineskip}
\newcommand{\setlinespacing}[1]%
{\setlength{\baselineskip}{#1 \defbaselineskip}}

\date{}
\begin{document}
\title{The spectrum of the Corona of Hypergraphs}
\author{Liya Jess Kurian$^1$\footnote{liyajess@gmail.com},  Chithra A. V$^1$\footnote{chithra@nitc.ac.in}
 \\ \small 
 1 Department of Mathematics, National Institute of Technology Calicut,\\\small
 Calicut-673 601, Kerala, India\\ \small	}
\maketitle
\thispagestyle{empty}
 \begin{abstract}
  The corona of hypergraphs is an extension of the corona operation applied to graphs.  The corona $G_0^* \odot_1^n G_1^*$ of two hypergraphs is obtained by taking $n$ copies of $G_1^*$ (where $n$ is the order of  $G_0^*$) and by joining the $i$-th vertex of $G_0^*$ with the $i$-th copy of $G_1^*$. In this paper, we estimate the complete spectrum(adjacency and Seidel) and eigenvectors of the corona $G_0^* \odot_1^n G_1^*$ of two hypergraphs when $G_1^*$ is regular. Additionally, we define the corona hypergraph $G_0^{*(m)}=G_0^{*(m-1)} \odot_1^n G_0^*$ and determined its adjacency spectrum. Also, we extend the definition coronal of the adjacency matrix. Moreover, we estimate the characteristic polynomial of Seidel matrix of the generalised corona of hypergraphs. Applying these results, we obtain infinitely many non-regular non-isomorphic adjacency and Seidel cospectral hypergraphs.
\end{abstract}
\section{Introduction}
Hypergraph is the generalisation of a graph in which hyperedge can connect more than 2 vertices at a time. In this article, we study hypergraphs that are simple and finite. More precisely, a hypergraph $G^*=(V,E)$ of order $n$ with a vertex set $V(G^*)=\{v_1,v_2,v_3,\ldots,v_n\}$ and a collection $E(G^*)=\{e_1,e_2,e_3,\ldots,e_m\}$ of hyperedges such that each hyperedge is a  non-empty subset $(|e_i| \geq 2)$ of the vertex set $V$. The degree $d(v)$ of a vertex $v\in V$ is defined as the number of hyperedges that include the vertex $v$. A hypergraph $G^*$ is said to be a $k$-uniform hypergraph \textnormal{\cite{Cooper2012,Kumar2017}} if each of its hyperedges contains exactly $k$  vertices, where $k\geq 2$. A hypergraph with $d(v_i)=r$ for all $v_i\in V$ is called a $r$-regular hypergraph. A hypergraph is said to be a $(k,r)$-regular hypergraph if it is both  $k$-uniform  and $r$-regular. Hypergraphs that are 2-uniform and 2-regular are referred to as cycles in graph theory. The properties of $(k,r)$ regular hypergraph are studied in \cite{Kumar2017}. A $k$-uniform hypergraph of order $n$ is said to be a complete $k$-uniform hypergraph $K_n^k$ if every possible $k$ subset of the vertex set forms a hyperedge \cite{Berge1973}. Let $G_1^*$  and $G_2^*$ be two $k$-uniform hypergraphs with vertex sets $V_1(G_1^*)$ and $V_2(G_2^*)$ and edge sets $E_1(G_1^*)$ and $E_2(G_2^*)$, respectively. Then the join of two $k$-uniform hypergraphs\cite{Sarkar2020} $G_1^*$ and $G_2^*$ is also a $k$-uniform hypergraph $G^*=G_1^*\oplus G^*_2$  with vertex set $V(G^*)=V_1(G_1^*)\cup V_2(G_2^*)$ and edge set $E(G^*)=E_0\cup \bigcup\limits_{i=1}^{2}E_i(G_i^*)$, where $E_0$ contains all possible $k$-subsets of $V$ satisfying $e \cap V_i\neq \phi,\, \forall i=1,2, \,e \in E_0$. Let $G_1^*[V_1']$ be the subhypergraph induced by the vertex set $V_1'\subset V_1(G_1^*)$. Then $G^*_1[V_1']\circledast G^*_2$ is a hypergraph with vertex set $V_1(G_1^*)\cup V_2(G_2^*)$ and edge set contains the hyperedges of $G_1^*$ together with hyperedges of $G_1^*[V_1']\oplus G_2^*$.\\

Research on the spectrum of hypergraphs has gained significant attention. In 2012, Cooper and Dutle \cite{Cooper2012} suggested tensors (hypermatrices) to analyse the spectral characteristics of hypergraphs, which gives analogous results to spectral graph theory. Here, we consider adjacency and Seidel matrices associated with the hypergraph. The adjacency matrix $A(G^*)=(a_{ij})$ of $G^*$ \cite{Lin2017} is an $ n\times n $  matrix whose rows and columns are indexed by the vertices of $G^*$ and for all $ v_i,v_j \in V, $
 \begin{equation*}
a_{ij} =\left\{
   \begin{array}{ll}
      \mid \{e_k  \in E(G^*): \{v_i,v_j\} \subset  e_k\}\mid    &  \mbox{, } v_i \neq v_j, k=1,2,3,...,m   \\
       0  & \mbox{, } v_i=v_j
   \end{array}.
   \right.
\end{equation*}
Let $J_n$ and $I_n$ denote the all one and identity matrix of order $n$ and $ J_{k,n}$ denote all one matrix of order $k\times n$, respectively. Then, the Seidel matrix $S(G^*)$ of a hypergraph $G^*$ is defined as $S(G^*)=J_n-I_n-2A(G^*)$ \cite{Zakiyyah2022}. For any square matrix $M$ we can find a scalar $\lambda$ such that $M\mathbf{x}=\lambda \mathbf{x} $ where $\mathbf{x}$ is the non-zero eigenvector corresponding to the eigenvalue $\lambda$. 
The collection of all the eigenvalues of $A(G^*)$ and $S(G^*)$ together with their multiplicities is known as the adjacency spectrum and the Seidel spectrum of $G^*$, respectively.  Let $\lambda_{1},\lambda_{2},\lambda_{3},\ldots,\lambda_{d},$ be the eigenvalues of an adjacency matrix $A(G^*)$ of hypergraph $G^*$ with multiplicities $m_{1},m_{2},m_{3},...,m_{d}$. Then, the adjacency spectrum $\sigma_{A}(G^*)$ of $G^*$ is denoted as,
$$ \sigma_{A}(G^*)=\begin{pmatrix}
               \lambda_{1} & \lambda_{2} & \lambda_{3} & \cdots & \lambda_{d}\\
               m_{1} & m_{2} & m_{3} & \cdots  & m_{d}
             \end{pmatrix}.$$
The adjacency spectral radius(Seidel spectral radius) of $G^*$, denoted by $\rho\left(A({G^*})\right)\,(  \rho\left(S({G^*})\right))$ is the largest eigenvalue of $A(G^*)\,( S(G^*))$. Two hypergraphs $G^*$ and $H^*$ are Seidel-cospectral(adjacency-cospectral) if they have the same Seidel spectrum(adjacency spectrum). For a $k$-uniform hypergraph $G^*=(V,E)$, $v\in V$ is said to be the neighbour of some $(k-1)$ subset of $V$, say $U$ if $\{v\}\cup U$ is a hyperedge of $G^*$. The neighbourhood of $U$ is denoted by $N(U)$.

In \cite{Frucht1970}, Frucht and Harary proposed the construction of corona graphs. Subsequently, the analysis of the spectra of the corona product of two graphs has received huge attention. In \cite{Barik2007}, Barik et al., have provided the complete spectra of the corona product of two graphs as a means of determining the spectrum of larger graphs in terms of simple graphs. By introducing the concept of coronal McLeman \cite{McLeman2011} formulated a method to establish the spectrum of the corona product of two hypergraphs. Later, Sharma\cite{Sharma2017} determined the structural properties of the corona graph and its spectrum. The study on the corona product of graphs was recently extended to hypergraphs \cite{Sarkar2020}. In connection with the spectrum of hypergraphs, it is of interest to identify a family of non-isomorphic cospectral hypergraphs. Obviously, isomorphic graphs are cospectral. There are several methods to find a cospectral family of hypergraphs. In this paper, we use signed M-GM switching to construct a family of non-isomorphic Seidel cospectral hypergraphs. Throughout this paper, we use the notation $j\in\mathbb{Z}^{[a,b]}$ to denote $j$ that takes all the integer values satisfying the condition $a\leq j\leq b$.

This paper focuses on the study of the corona product of $k$-uniform hypergraphs. In Section 2, we give fundamental definitions and results that will be used later. In Section 3, we determine the adjacency and seidel spectrum of generalised corona of $k$-uniform hypergraphs. In Section 4, we estimate the adjacency and Seidel spectrum of corona product of two $k$-uniform hypergraphs. Also, we calculate the coronal of $(k,r)$-regular hypergraph and find an infinite pair of adjacency cospectral hypergraphs. In Section 5, we compute the adjacency spectrum of the corona hypergraph. In Section 6, we present a new method for the construction of Seidel-cospectral hypergraphs. Using this we get infinitely many non-isomorphic Seidel-cospectral $k$-uniform hypergraphs.
\section{Preliminaries}

 \begin{lemma}\label{detblock}\textnormal{\cite{Das2018}}
Let $M_{11}, M_{12}, M_{21}$, and $M_{22}$ be matrices with $M_{22}$ invertible. Let 
             $$ M= \begin{bmatrix}
                 M_{11} & M_{12}\\
                M_{21} & M_{22}
                  \end{bmatrix}$$
Then, $det(M) = det(M_{22})det(M_{11}-M_{12}M_{22}^{-1}M_{21}).$
          \end{lemma}
 \begin{lemma}\textnormal{\cite{jahfar2020}}\label{detij}
              For any two real numbers $r$ and $s$,
              $$(rI_n-sJ_n)^{-1}=\frac{1}{r}I_n+\frac{s}{r(r-ns)}J_n.$$
\end{lemma}
\begin{lemma}\textnormal{\cite{Ding2007}}\label{matrixdet}
Let $M\in\mathbb{R}^{n\times n}$ be an invertible  matrix, $U$ and $W$ are $n \times 1$ matrices. Then
$$\det(UW^{T}+M)=(1+W^{T}M^{-1}U)\det(M).$$
\end{lemma}
\begin{lemma}\cite{Sharma2017}\label{noofeig}
  Let $n,\,m\in \mathbb{R}$. Then
  $$\sum_{j=0}^{m-1} 2^jn(n-1)(n+1)^{m-j-1}+2^mn=n(n+1)^m.$$
\end{lemma}
\begin{definition}\cite{Sarkar2020}(Generalized corona of hypergraphs)\label{defn:GenCoro}
    The corona $G^*=G_0^* \odot_p^t G_i^*, ~i\in\mathbb{Z}^{[1,t]} $ of $k$-uniform hypergraphs $G_0^*$ and $G_i^*$ (where $G_0^*$ and $G_i^*$  has $n$ and $n_i$ vertices respectively) with partition $P=\{U_1,U_2,\ldots,U_t\}(\text{with}~ |U_i|=p,\,p\in \mathbb{N})$ of  the vertex set $V_0$ of hypergraph $G_0^*$ is defined as a $k$-uniform hypergraph  $\bigcup\limits_{i=1}^t \Bigl(\bigcup\limits_{j=1}^p G_0^*[U_i]\circledast G_i^{*(j)}\Bigr)  $ obtained by taking $p$ copies of $G_i^*$ say $G_i^{*(j)}$, where $j\in\mathbb{Z}^{[1,p]}$. 
\end{definition}
\section{Generalised corona of hypergraphs}

 In this section, we determine the spectral properties of the generalised corona of $k$-uniform hypergraphs. In \cite{Sarkar2020}, Sarkar and Banerjee provided the spectrum of the generalised corona of $k$-uniform hypergraphs by considering another definition of the adjacency matrix. We modify the spectrum by taking into account the definition stated in \cite{Lin2017}.

\subsection{Adjacency spectrum of corona of hypergraphs}
Let $G_0^*$ be a $k$-uniform hypergraph with a vertex set $V_0$ and partition $P=\{U_1,U_2,\ldots,U_t\}$ such that $|U_i|=p,\,p\in \mathbb{N} $. Also, $G_i^*$ be the $(k,r)$- regular hypergraphs with vertex set $V_i$ and edge set $E_i$, where $i\in\mathbb{Z}^{[1,t]}$. Let $|V_i|=m$, then the corresponding adjacency matrix of $G^*$ given by

\begin{equation*}
    A(G^*)=\begin{bmatrix}
        X & H\\
        H^T & Y
    \end{bmatrix},
\end{equation*}
where $X$, $H$ and $Y$ denotes the adjacency relation between the vertices of $V_0$, vertices of $V_0$ and $V_i$, and vertices of $V_i$ respectively.
\begin{align*}
    X&=A(G_0^*)+aI_t\otimes (J_p-I_p), \\
   H&=bI_t\otimes J_{p, pm},\\
   Y&=diag(I_p\otimes Y_1, I_p\otimes Y_2,I_p\otimes Y_3,\ldots,I_p\otimes Y_t),
\end{align*}
where $a=\begin{cases}
			p\Bigl(\binom{p+m-2}{k-2}-\binom{p-2}{k-2}\Bigr) & \text{if $p\in \mathbb{Z}^{[2,n]}$ }\\
            0 & \text{if $p=1$}
		 \end{cases},\,b=\binom{p+m-2}{k-2},\,Y_i=A(G_i^*)+c(J_m-I_m)$ and $c=\binom{p+m-2}{k-2}-\binom{m-2}{k-2}$. 
\begin{thm}\label{thm3.1adj}
   The characteristic polynomial of  $G^*=G_0^* \odot_p^t G_i^*,\,i\in \mathbb{Z}^{[1,t]}$ is given by,
   $$ P_{A(G^*)}(\lambda)=\left(\prod \limits_{i=1}^t \det(Y_i-\lambda I_m)\right)^p \det(A(G_0^*)+I_t\otimes \bigl((a-\frac{b^2pm}{r(k-1)+c(m-1)-\lambda})J_p-(a+\lambda)I_p\bigr),$$
   where $a=\begin{cases}
			p\Bigl(\binom{p+m-2}{k-2}-\binom{p-2}{k-2}\Bigr) & \text{if $p\in \mathbb{Z}^{[2,n]}$ }\\
            0 & \text{if $p=1$}
		 \end{cases},\,b=\binom{p+m-2}{k-2},\,Y_i=A(G_i^*)+c(J_m-I_m)$ and $c=\binom{p+m-2}{k-2}-\binom{m-2}{k-2}$.
\end{thm}
\begin{proof}
    The characteristic polynomial of $G^*$ is given by,
    \begin{equation*}
    P_{A(G^*)}(\lambda)=\det\begin{pmatrix}
        A(G_0^*)+aI_t\otimes (J_p-I_p)-\lambda I_n & b I_t\otimes J_{p, pm}\\
        b I_t\otimes J_{pm, p} & Y-\lambda I_{ptk}
    \end{pmatrix}.
\end{equation*}
From Lemma \ref{detblock}, we get
\begin{equation}\label{eqn1.1}
     P_{A(G^*)}(\lambda)=\det(Y-\lambda I_{ptk})\det(A(G_0^*)+aI_t\otimes (J_p-I_p)-\lambda I_n- b^2 I_t\otimes J_{p, pm}(Y-\lambda I_{ptk})^{-1} I_t\otimes J_{pm, p}). 
\end{equation}
Since row sum of the adjacency matrix of a $(k,r)-$regular hypergraph is $r(k-1)$, we get
\begin{align*}
    (Y-\lambda I_{ptk})(I_t\otimes J_{pm, p})&=(r(k-1)+c(m-1)-\lambda)(I_t\otimes J_{pm, p}).
\end{align*}
Then,
\begin{equation}\label{eqn1.2}
(I_t\otimes J_{pm, p})=\frac{1}{(r(k-1)+c(m-1)-\lambda)}(Y-\lambda I_{ptk})(I_t\otimes J_{pm, p}).
\end{equation}
Also,
\begin{equation}\label{eqn1.3}
(I_t\otimes J_{pm, p})(I_t\otimes J_{pm, p})=pm I_t\otimes J_p.
\end{equation}
Hence the theorem follows from (\ref{eqn1.1}),(\ref{eqn1.2}) and (\ref{eqn1.3}).
\end{proof}
Using Theorem \ref{thm3.1adj}, we obtain the characteristic polynomial of $G^*$ in terms of the eigenvalues of the factor hypergraphs $G^*_0$ and $G_i^*, i\in\mathbb{Z}^{[1,n]}$ when $p=1$ and $t=n$.
\begin{corlly}
Let $G_0^*$ be a $k$-uniform hypergraph with eigenvalues $\lambda_1'\geq \lambda_2'\geq \lambda_3'\geq \ldots \geq \lambda_n'$ and $G_i^*,i\in\mathbb{Z}^{[1,n]}$ be $(k,r)$-regular hypergraphs with eigenvalues $\lambda_1^{(i)}=r(k-1)\geq \lambda_2^{(i)}\geq \lambda_3^{(i)}\geq \ldots \geq \lambda_m^{(i)}$. Then the characteristic polynomial of  $G^*=G_0^* \odot_1^n G_i^*$ is given by,
\begin{align*}
            P_{A(G^*)}(\lambda)=(r(k-1)-c(m-1)-\lambda)^n\prod \limits_{i=1}^n\prod \limits_{j=2}^m  (\lambda_j^{(i)}-c-\lambda)\prod \limits_{i=1}^n( \lambda_{i}'+\frac{\lambda^2-(r(k-1)+c(m-1))\lambda-b^2m}{r(k-1)+c(m-1)-\lambda}).
       \end{align*}
\end{corlly}
\begin{proof}
     By substituting $t=n$ and $p=1$ in Theorem \ref{thm3.1adj}, we have
       $$ P_{A(G^*)}(\lambda)=\left(\prod \limits_{i=1}^n \det(Y_i-\lambda I_m)\right) \det\biggl(A(G_0^*)+\bigl(\frac{\lambda^2-(r(k-1)+c(m-1))\lambda-b^2m}{r(k-1)+c(m-1)-\lambda}\bigr)I_n \biggr).$$
       Hence the result.
      
\end{proof}
\subsection{Seidel Spectrum of the corona of hypergraphs}
 Let $ G_0^*$ be a $k-$uniform hypergraph of order $n$ with the partition $P=\{U_1,U_2,U_3,\ldots, U_t\} $, where $|U_i|=p,\, p\in \mathbb{N}$ and $G_i^*,\,i\in\mathbb{Z}^{[1,t]}$ be a $k$-uniform hypergraph of order $n_i$. Then the Seidel matrix of $G^*=G_0^*\odot_p^tG_i^*$ is given by,
  $$ S(G^*)=\begin{bmatrix}
	X_S & H_S\\
	H_S^T & Y_S
\end{bmatrix}$$

such that,
\begin{align*}
	X_S&=S(G_0^*)-2aI_t\otimes (J_p-I_p),\\
	H_S&=(J_t-2bI_t )\otimes J_{p,pn_i},\\
	Y_S&=J_{pn_it}+diag(I_p\otimes(Y_{S_1}-J_{n_i}),\ldots,I_t\otimes(Y_{S_t}-J_{n_i})),
\end{align*}
where $Y_{S_i}=S(G_i^*)-2c(J_{n_i}-I_{n_i}) $.
\begin{thm}\label{thm2.4}
	Let $ G_0^*$ be a $k-$uniform hypergraph of order $n$ with the partition $P=\{U_1,U_2,U_3,\ldots, U_t\} $, where $|U_i|=p,\, p\in \mathbb{N}$ and $G_i^*,\,i\in\mathbb{Z}^{[1,t]}$ be a $(k,r)$-regular hypergraph of order $m$. Then the characteristic polynomial of Seidel matrix of $G^*=G_0^* \odot_p^t G_i^*$ is
 \begin{equation*}
    \begin{split}
        P_{S(G^*)}(\mu)=(1+\frac{pmt}{h})\biggl( \prod_{i=1}^t \det (Y_{S_i}-J_m-\mu I_m)\biggr)^p \det\biggl(S(G_0^*)+I_t\otimes \Bigl(-(2a+\frac{4pmb^2}{h})J_p-I_p)\\+(2a-\mu) I_p\Bigr)-\frac{pm}{h} \Bigl((t-4b)-\frac{pm(t-2b)^2}{h+pmt} \Bigr) J_{n}\biggr),
    \end{split} 
 \end{equation*}
 where $a=\begin{cases}
			p\Bigl(\binom{p+m-2}{k-2}-\binom{p-2}{k-2}\Bigr) & \text{if $p\in \mathbb{Z}^{[2,n]}$ }\\
            0 & \text{if $p=1$}
		 \end{cases},\,b=\binom{p+m-2}{k-2},\,Y_{S_i}=S(G_i^*)-2c(J_m-I_m)$ , $c=\binom{p+m-2}{k-2}-\binom{m-2}{k-2}$ and $h=-(1+2r(k-1)+2c(m-1)+\mu)$.
\end{thm}
\begin{proof}
	The characteristic polynomial of the Seidel matrix of $G^*$ is given by
	$$ P_{S(G^*)}(\mu)=\det\begin{pmatrix}
		S(G_0^*)-2aI_t\otimes (J_p-I_p)-\mu I_n & H_S\\
		H_S^{T} & Y_S-\mu I_{pmt}
	\end{pmatrix},$$
	where $Y_S-\mu I_{pmt}= J_{pmt}+diag(I_p\otimes(Y_{S_1}-J_m-\mu I_m),\ldots,I_p\otimes(Y_{S_t}-J_m-\mu I_m))$ and $H_S=(J_t-2bI_t )\otimes J_{p,pm}.$\\
	From Lemma \ref{detblock}, we get 
		\begin{equation}\label{eqn4a}
		    P_{S(G^*)}(\mu)=\det(Y_S-\mu I)\det(S(G_0^*)-2aI_t\otimes (J_p-I_p)-\mu I_n-H_S(Y_S-\mu I_{pmt})^{-1}H_S^T).
		\end{equation}
	Since $G_i^*$'s are $r$-regular $m$-uniform hypergraphs, we have
	$$ S(G_i^*)J_{m,p}=(m-1-2r(k-1))J_{m,p}.$$
	On simplification, we get 
	\begin{align}\label{eqn4}
 \begin{split}
 \begin{rcases}
  diag(I_p\otimes(Y_{S_1}-J_m-\mu I_m),\ldots,I_p\,\otimes&(Y_{S_t}-J_m-\mu I_m))(J_t\otimes J_{m,p})\\
                                                                            &=-(1+2r(k-1)+2c(m-1)+\mu)(J_t\otimes J_{m,p}),\\
		diag(I_p\otimes(Y_{S_1}-J_m-\mu I_m),\ldots,I_p\,\otimes&(Y_{S_t}-J_m-\mu I_m))(I_t\otimes J_{m,p})\\
                                                            &=-(1+2r(k-1)+2c(m-1)+\mu)(I_t\otimes J_{m,p}).    
 \end{rcases}
 \end{split}
	\end{align}
	Now consider,
 \begin{equation*}
     (Y_S-\mu I_{pmt})H_S^T=\bigl( J_{pmt}+diag(I_p\otimes(Y_{S_1}-J_m-\mu I_m),\ldots,I_p\otimes(Y_{S_t}-J_m-\mu I_m))\bigr)\bigl((J_t-2bI_t )\otimes J_{pm,p}\bigr).
 \end{equation*}
 From (\ref{eqn4}), we get
 \begin{align*}
      (Y_S-\mu I_{pmt})H_S^T&=J_{pmt}\bigl((J_t-2bI_t )\otimes J_{pm,p}\bigr)-(1+2r(k-1)+2c(m-1)+\mu)\bigl((J_t-2bI_t )\otimes J_{pm,p}\bigr)\\
      &=\bigl(J_{pmt}+hI_{pmt}\bigr)H_S^T,
 \end{align*}

where $h=-(1+2r(k-1)+2c(m-1)+\mu)$. Therefore,
 \begin{equation*}
    H_S^T=\bigl(J_{pmt}+hI_{pmt}\bigr)^{-1}(Y_S-\mu I_{pmt})H_S^T. 
 \end{equation*}
 By Lemma \ref{detij}, we get $\bigl(hI_{pmt}+J_{pmt}\bigr)^{-1}=\frac{1}{h}I_{pmt}-\frac{1}{h(h+pmt)}J_{pmt}$. Hence,
\begin{equation}\label{eqn5}
    H_S^T=\bigl(\frac{1}{h}I_{pmt}-\frac{1}{h(h+pmt)}J_{pmt}\bigr)(Y_S-\mu I_{pmt})H_S^T.
\end{equation}
Again,
\begin{align}\label{eqn6}
    H_SH_S^T&=\Bigl((J_t-2bI_t )\otimes J_{p,pm}\Bigr)\Bigl((J_t-2bI_t )\otimes J_{pm,p}\Bigr)\notag \\
        &=\bigl( (t-4b)J_t+4b^2I_t \bigr) \otimes pmJ_{p}.
\end{align}
Also, 
\begin{align}\label{eqn7}
    H_SJ_{pmt}H_S^T&=\Bigl((J_t-2bI_t )\otimes J_{p,pm}\Bigr)J_{pmt}\Bigl((J_t-2bI_t )\otimes J_{pm,p}\Bigr)\notag \\
    &=\Bigl(J_{pt,pmt}J_{pmt}-2b(I_t\otimes J_{p,pm})J_{pmt} \Bigr)\Bigl((J_t-2bI_t )\otimes J_{pm,p}\Bigr)\notag \\
    &=( pmt-2bpm)J_{pt,pmt} \Bigl((J_t-2bI_t )\otimes J_{pm,p}\Bigr)\notag \\
    &=( pmt-2bpm)^2 J_{pt}=p^2m^2(t-2b)^2 J_{pt}.
\end{align}
From (\ref{eqn5}), we have
\begin{align*}
    H_S(Y_S-\mu I_{pmt})^{-1}H_S^T&=H_S(Y_S-\mu I_{pmt})^{-1}\bigl(\frac{1}{h}I_{pmt}-\frac{1}{h(h+pmt)}J_{pmt}\bigr)(Y_S-\mu I_{pmt})H_S^T\\
    &=\frac{1}{h}H_SH_S^T-\frac{1}{h(h+pmt)}H_S(Y_S-\mu I_{pmt})^{-1}J_{pmt}(Y_S-\mu I_{pmt})H_S^T.
\end{align*}
Since $Y_S-\mu I_{pmt}$ commute with $J_{pmt}$, we get
\begin{align*}
    H_S(Y_S-\mu I_{pmt})^{-1}H_S^T&=\frac{1}{h}H_SH_S^T-\frac{1}{h(h+pmt)}H_SJ_{pmt}H_S^T.
\end{align*}
From (\ref{eqn6}) and (\ref{eqn7}), we obtain

\begin{equation}\label{eqn9}
     H_S(Y_S-\mu I_{pmt})^{-1}H_S^T=\frac{pm}{h}\biggl( \Bigl((t-4b)-\frac{pm(t-2b)^2}{h+pmt} \Bigr) J_{pt}+ 4b^2(I_t \otimes J_{p})\biggr).
\end{equation}
Take $D=diag(I_p\otimes(Y_{S_1}-J_m-\mu I_m),\ldots,I_p\otimes(Y_{S_t}-J_m-\mu I_m))$,  Then 
$DJ_{pmt,1}=hJ_{pmt,1}$. Also,
\begin{equation}\label{eqn10}
    J_{pmt,1}=\frac{1}{h}DJ_{pmt,1}.
\end{equation}
Now consider,
$$\det(Y_S-\mu I)=\det(J_{pmt}+D)=\det(J_{pmt,1}J_{1,pmt}+D).$$
By Lemmma \ref{matrixdet}, we have
$$\det(Y_S-\mu I)=(1+J_{1,pmt}D^{-1}J_{pmt,1})\det(D).$$
Using (\ref{eqn9})  and (\ref{eqn10}) in  (\ref{eqn4a}), we get
\begin{equation*}
    \begin{split}
         P_{S(G^*)}(\mu)=(1+\frac{pmt}{h})\biggl( \prod_{i=1}^t \det (Y_{S_i}-J_m-\mu I_m)\biggr)^p \det\biggl(S(G_0^*)+I_t\otimes \Bigl(-(2a+\frac{4pmb^2}{h})J_p-I_p)\\+(2a-\mu) I_p\Bigr)-\frac{pm}{h} \Bigl((t-4b)-\frac{pm(t-2b)^2}{h+pmt} \Bigr) J_{n}\biggr).
    \end{split}
\end{equation*}

\end{proof}
Next, we determine the characteristic polynomial of $G^*$ in terms of the factor hypergraphs when $p=1$ and $t=n$.
\begin{thm}\label{S-cosepec_1}
  	Let $ G_0^*$ be a $(k,r_0)$-regular hypergraph of order $n$ with partition $P=\{U_1,U_2,U_3,\ldots, U_t\} $, where $|U_i|=p, p\in \mathbb{N}$ and $G_i^*, i \in \mathbb{Z}^{[1,t]}$ be a $(k,r)$-regular hypergraph of order $m$. Let $\mu_1^{(0)}=n-1$ $-2r_0(k-1),\,$  $ \mu_2^{(0)}\geq \mu_3^{(0)}\geq \ldots \geq \mu_n^{(0)}$ be the Seidel eigenvalues of $G_0^*$ and $\mu_1^{(i)}=m-1-2r(k-1),\,$  $ \mu_2^{(i)}\geq \mu_3^{(i)}\geq \ldots \geq \mu_m^{(i)}$ be the Seidel eigenvalues of $G_i^*$. Then the characteristic polynomial of the Seidel matrix of $G^*=G_0^* \odot_1^n G_i^*$ is
 \begin{equation*}
    \begin{split}
        P_{S(G^*)}(\mu)=\bigl((\mu_1^{(0)}-\mu)(h+mn)+4bm(n-b)-mn^2\bigr)\Bigl(\prod_{i=2}^n(\mu_i^{(0)}h-4mb^2-\mu h)\Bigr)\Bigl(\prod_{j=1}^n \prod_{i=2}^m (\mu_i^{(j)}+2c-\mu)\Bigr)
    \end{split} 
 \end{equation*}
 where $b=\binom{p+m-2}{k-2},\,c=\binom{p+m-2}{k-2}-\binom{m-2}{k-2}$ and $h=-(1+2r(k-1)+2c(m-1)+\mu)$.  
\end{thm}
\begin{proof}
    By substituting $t=n$ and $p=1$ in Theorem \ref{thm2.4}, we obtain
    \begin{equation}\label{eqn11*}
    \begin{split}
        P_{S(G^*)}(\mu)
        =\bigl(\frac{h+mn}{h}\bigr)\biggl( \prod_{i=1}^n \det (Y_{S_i}-J_m-\mu I_m)\biggr) \det\biggl(S(G_0^*)+\Bigl(\frac{-4mb^2-\mu h}{h}\Bigr)I_n-\frac{m}{h} \Bigl(\frac{nh-4bh-4b^2m}{h+mn} \Bigr) J_{n}\biggr).
     \end{split} 
 \end{equation}
 Consider,
 \begin{align}\label{eqn12*}
     \det\biggl(S(G_0^*)&+\Bigl(\frac{-4mb^2-\mu h}{h}\Bigr)I_n-\frac{m}{h} \Bigl(\frac{nh-4bh-4b^2m}{h+mn} \Bigr) J_{n}\biggr)\notag\\
     &=\Bigl(\mu_1^{(0)}-\frac{mn(nh-4bh-4b^2m)}{h(h+mn)}-\frac{4mb^2+\mu h}{h}\Bigr)\Bigl(\prod_{i=2}^n(\mu_i^{(0)}-\frac{4mb^2+\mu h}{h})\Bigr)\notag\\
     &=\Bigl(\frac{\mu_1^{(0)}h(h+mn)-mn(nh-4bh-4b^2m)-(4mb^2+\mu h)(h+mn)}{h^n(h+mn)}\Bigr)\Bigl(\prod_{i=2}^n(\mu_i^{(0)}h-4mb^2+\mu h)\Bigr)\notag\\
     &=\bigl(\frac{(\mu_1^{(0)}-\mu)(h+mn)+4bm(n-b)-mn^2}{h^n(h+mn)}\bigr)\Bigl(\prod_{i=2}^n(\mu_i^{(0)}h-4mb^2+\mu h)\Bigr).
 \end{align}
 Also, we have
 \begin{align}\label{eqn13*}
     \prod_{i=1}^n \det(Y_{S_i}-J_m-\mu I_m)&=\prod_{i=1}^n \det(S(G_i^*)-(2c+1)J_m+(2c-\mu) I_m)\notag\\
     &=(m-1-2r(k-1)-(2c+1)m+2c-\mu)^n \prod_{j=1}^n \prod_{i=2}^m (\mu_i^{(j)}+2c-\mu)\notag\\
     &=h^n \prod_{i=2}^m (\mu_i^{(j)}+2c-\mu).
 \end{align}
 Substituting (\ref{eqn12*}),(\ref{eqn13*}) in (\ref{eqn11*}), we get the desired result.
\end{proof}
\section{Spectrum of the corona of two hypergraphs}

The corona product of two $k$-uniform hypergraphs $G^*=G_0^* \odot_1^n G_1^*$, is the $k$-uniform hypergraph obtained by taking $p=1$, $k=n$, and $G_i^*=G_1^*$ for $i\in\mathbb{Z}^{[1,n]}$ in Definition \ref{defn:GenCoro}. By proper labelling of vertices, we obtain the adjacency matrix and Seidel matrix of $G^*$ as follows:
\begin{equation*}
    A(G^*)=\begin{bmatrix}
        A(G_0^*) & \binom{m-1}{k-2}J_{1,m}\otimes I_n\\
        \binom{m-1}{k-2}J_{m,1}\otimes I_n & A(G_1^*)\otimes I_n
    \end{bmatrix},
\end{equation*}
\begin{equation*}
    S(G^*)=\begin{bmatrix}
        S(G_0^*)&\displaystyle J_{1,m}\otimes\left(J_n-2\binom{m-1}{k-2}I_n\right) \\
       \displaystyle J_{m,1}\otimes\left(J_n-2\binom{m-1}{k-2}I_n\right) &J_m\otimes (J_n-I_n)+S(G_1^*)\otimes I_n
    \end{bmatrix}.
\end{equation*}
The following theorem describes the adjacency spectrum and spectral radius of corona product of two $k$-uniform hypergraphs $G_0^*$ and $G_1^*$, where $G_1^*$ is $r$-regular hypergraph.
\begin{theorem}\label{corona2hygrph}
    Let $G_0^*$ be any $k$-uniform hypergraph of order $n$ with eigenvalues $\lambda_1^{(0)}\geq \lambda_2^{(0)}\geq \lambda_3^{(0)}\geq \ldots \geq \lambda_n^{(0)}$and $G_1^*$ be a $(k,r)$-regular hypergraph of order $m$ with eigenvalues $\lambda_1^{(1)}=r(k-1)\geq \lambda_2^{(1)}\geq \lambda_3^{(1)}\geq \ldots \geq \lambda_m^{(1)}$, then for $i\in\mathbb{Z}^{[1,n]}, j\in\mathbb{Z}^{[2,m]}$ the spectrum of $G^*=G_0^* \odot_1^n G_1^*$ is given by
     $$\sigma_A({G^*})=
    \begin{pmatrix}
     \frac{\lambda_i^{(0)}+r(k-1)\pm \sqrt{\left(\lambda_i^{(0)}-r(k-1)\right)^2+4m\mbinom{m-1}{k-2}^2}}{2}&\lambda_j^{(1)} \\
                 1 &n
    \end{pmatrix}$$
    and $\rho_A(G^*)=\displaystyle  \frac{\rho(A(G^*_0))+r(k-1)+ \sqrt{\left(\rho(A(G^*_0))-r(k-1)\right)^2+4m\mbinom{m-1}{k-2}^2}}{2}.$
\end{theorem}
\begin{proof}
For $i\in \mathbb{Z}^{[1,n]}$, let $\mathbfit{x}_i$ be an eigenvector of $A(G_0^*)$ corresponding to the eigenvalue $ \lambda_i^{(0)}$. We have to determine $\lambda_i$'s and $t_i$'s such that $\mathbf{x}_i=\begin{bmatrix}
    \mathbfit{x}_i & t_iJ_{m,1}\otimes\mathbfit{x}_i
                       \end{bmatrix}^T$
is an eigenvector corresponding to the eigenvalue $\lambda_i$ of $A(G^*)$. Since $A(G^*)\mathbf{x}_i=\lambda_i\mathbf{x}_i$, we have the following equations
\begin{align*}
    \lambda_i^{(0)} \mathbfit{x}_i+t_i\mbinom{m-1}{k-2}(J_{1,m}\otimes I_n)(J_{m,1}\otimes \mathbfit{x}_i)&=\lambda_i \mathbfit{x}_i,\\
    \mbinom{m-1}{k-2}(J_{m,1}\otimes I_n) \mathbfit{x}_i+t_i(A(G_1^*)\otimes I_n)(J_{m,1}\otimes \mathbfit{x}_i)&=\lambda_it(J_{m,1}\otimes \mathbfit{x}_i).
\end{align*}
On simplification, we get
\begin{align*}
    \lambda_i^{(0)}+t_i m \mbinom{m-1}{k-2} =\lambda_i,\;\;\; \;\mbinom{m-1}{k-2}+t_ir(k-1)=\lambda_i t_i.
\end{align*}
  
After solving these equations, we obtain
$$
 \lambda_i=\frac{\lambda_i^{(0)}+r(k-1)\pm \sqrt{\left(\lambda_i^{(0)}-r(k-1)\right)^2+4m\mbinom{m-1}{k-2}^2}}{2},
$$
 and
$$t_i=\frac{-(\lambda_i^{(0)}-r(k-1))\pm \sqrt{\left(\lambda_i^{(0)}-r(k-1)\right)^2+4m\mbinom{m-1}{k-2}^2}}{2m\mbinom{m-1}{k-2}}=\frac{\binom{m-1}{k-2}}{\lambda_i-r(k-1)}. $$  
Hence, for $i\in\mathbb{Z}^{[1,n]} $ 
\begin{center}
   $ \Tilde{\lambda}_i=\frac{\lambda_i^{(0)}+r(k-1)+\sqrt{\left(\lambda_i^{(0)}-r(k-1)\right)^2+4m\mbinom{m-1}{k-2}^2}}{2},\; \Tilde{\Tilde{\lambda}}_i=\frac{\lambda_i^{(0)}+r(k-1)- \sqrt{\left(\lambda_i^{(0)}-r(k-1)\right)^2+4m\mbinom{m-1}{k-2}^2}}{2}$ 
\end{center}
are the eigenvalues of $A(G^*)$ corresponding to the eigenvectors
$$ \begin{bmatrix}
    \mathbfit{x}_i \\ \frac{\binom{m-1}{k-2}}{\Tilde{\lambda}_i-r(k-1)} J_{m,1}\otimes\mathbfit{x}_i
                       \end{bmatrix},\;\;
    \begin{bmatrix}
    \mathbfit{x}_i \\ \frac{\binom{m-1}{k-2}}{\Tilde{\Tilde{\lambda}}_i-r(k-1)} J_{m,1}\otimes\mathbfit{x}_i
                       \end{bmatrix}$$
respectively. Let $\mathbfit{y}_j$ be an eigenvector corresponding to the eigenvalue $\lambda_j^{(1)} $ of $A(G_1^*)$  and let $\mathbf{e}_i\in \mathbb{R}^n$ be a vector with 1 in the $i$-th coordinate and 0 elsewhere. Since $G_1^*$ is a $(k,r)$-regular hypergraph, for $j\in\mathbb{Z}^{[2,m]}$, $\mathbf{y}_j$'s are orthogonal to $\mathbfit{j}$. It is clear that, for $i\in\mathbb{Z}^{[1,n]}$
$$A(G^*)\begin{bmatrix}
    \mathbf{0}_{n,1} \\ \mathbfit{y}_j\otimes e_i
\end{bmatrix} = \lambda_j^{(1)} 
\begin{bmatrix}
    \mathbf{0}_{n,1} \\ \mathbfit{y}_j\otimes e_i
\end{bmatrix}. $$
For $2\leq j \leq m$,  $\lambda_j^{(1)}$ is an eigenvalue of $A(G^*)$ with multiplicity $n$. Note that $\lambda_j^{(1)}\leq r(k-1)$, hence $\rho(A(G^*))=\Tilde{\lambda}_1$. This completes the proof. 
\end{proof}
The following corollary gives the adjacency spectrum and spectral radius of $G^*$ when $G_1^*=K_m^m$. 
\begin{corollary}
     Let $G_0^*$ be any $m$-uniform hypergraph of order $n$ with eigenvalues $\lambda_1^{(0)}\geq \lambda_2^{(0)}\geq \lambda_3^{(0)}\geq \ldots \geq \lambda_n^{(0)}$and $G_1^*=K_m^m$ be an $m$-uniform complete hypergraph of order $m$, then for $i\in\mathbb{Z}^{[1,n]},\,j\in\mathbb{Z}^{[2,m]}$ the spectrum of $G^*=G_0^* \odot_1^n G_1^*$ is given by
     $$\sigma_A({G^*})=
    \begin{pmatrix}
     \frac{\lambda_i^{(0)}+m-1\pm \sqrt{\left(\lambda_i^{(0)}-(m-1)\right)^2+4m(m-1)^2}}{2}&-1 \\
                 1 &n(m-1)
    \end{pmatrix}$$
    and $\rho(A(G^*))=\displaystyle  \frac{\rho(A(G^*_0))+(m-1)+ \sqrt{\left(\rho(A(G^*_0))-(m-1)\right)^2+4m(m-1)}}{2}.$
    
\end{corollary}
\begin{example}
    Let $G_0^*=G_1^*=K_m^m$ be the complete $m$ uniform hypergraph of order $m$. Then spectrum of $G^*=K_m^m \odot_1^n K_m^m$ is
    $$\sigma_A({G^*})=
    \begin{pmatrix}
     (m-1)\left( 1\pm \sqrt{m} \right)&\displaystyle \frac{m-2\pm \sqrt{m^2+4m(m-1)^2}}{2} & -1\\
                 1 & m-1 & m(m-1)
    \end{pmatrix}$$
    and 
    $\rho(A(G^*))=(m-1)( 1 +\sqrt{m}).$
    \end{example}

The following theorem gives the Seidel spectrum of corona of two $k$-uniform hypergraphs.
\begin{theorem}\label{S-cospec2}
    Let $G_0^*$ be an $n$ order $(k,r_0)$-regular hypergraph with Seidel eigenvalues 
    $\mu_1^{(0)}\geq \mu_2^{(0)}\geq \mu_3^{(0)}\geq \ldots \geq \mu_n^{(0)}$and $G_1^*$ be $m$-order $(k,r_1)$-regular hypergraph  with  Seidel eigenvalues $\mu_1^{(1)}\geq \mu_2^{(1)}\geq \mu_3^{(1)}\geq \ldots \geq \mu_m^{(1)}$, then for $i\in \mathbb{Z}^{[1,n]},\,j\in \mathbb{Z}^{[1,m]},$ $\mu_i^{(0)}\neq n-1-2r_0(k-1)$ and $ \mu_j^{(1)}\neq m-1-2r_1(k-1)$ the Seidel spectrum of $G^*=G_0^* \odot_1^n G_1^*$ is given by
     $$\sigma_S({G^*})=
    \begin{pmatrix}
     \frac{\mu_i^{(0)}-1-2r_1(k-1)\pm \sqrt{\left(\mu_i^{(0)}+1+2r_1(k-1)\right)^2+16 m\mbinom{m-1}{k-2}^2}}{2}&\mu_j^{(1)} &\alpha_1 &\alpha_2\\
                 1 &n &1 &1
    \end{pmatrix},$$ where $\alpha_1$ and $ \alpha_2$ are the roots of the equation
     $\mu^2+\bigl(2-n(m+1)+2(r_0+r_1)(k-1)\bigr)\mu+\bigl(n-1-2r_0(k-1)\bigr)\bigl(mn-1-2r_1(k-1)\bigr)-m\bigl(n-2\mbinom{m-1}{k-2}\bigr)^2=0$.
\end{theorem}
\begin{proof}
      Let $\mathbfit{x}_i$ be an eigenvector associated with an eigenvalue $\mu_i^{(0)}(\mu_i^{(0)}\neq n-1-2r_0(k-1))$ and $i\in\mathbb{Z}^{[1,n]})$ of $S_{G_0}^*$. We have to determine $\mu_i$'s and $t_i$'s such that  $\mathbf{x}_i=\begin{bmatrix}
        \mathbfit{x}_i & t_i\,J_{m,1}\otimes \mathbfit{x}_i
    \end{bmatrix}^T$ is an eigenvector associated with the eigenvalue $\mu_i$ of $S_G^*$. Since   $S(G^*)\mu_i = \mu_i\mathbf{x}_i$, we have the following equations,
    \begin{align*}\label{eqn11a}
 \mu_i^{(0)}\mathbfit{x}_i+t_i m\Bigl(J_n-2\binom{m-1}{k-2}I_n\Bigr)\mathbfit{x}_i &=\mu_i\mathbfit{x}_i , \\
 \Bigl( J_{m,1}\otimes(J_n-2\binom{m-1}{k-2}I_n\Bigr)\mathbfit{x}_i+t_i\Bigl(mJ_{m,1}\otimes (J_n-I_n)\mathbfit{x}_i+S(G_1^*)J_{m,1}\otimes\mathbfit{x}_i\Bigr)&=\mu_it_iJ_{m,1}\otimes\mathbfit{x}_i .
	\end{align*}
 Also, row sum of $S(G_1^*)$ is $m-1-2r_1(k-1)$ and $\mathbfit{x}_i$ is orthogonal to $\mathbfit{j}$. Therefore, 
 \begin{align*}
     \mu_i&=\mu_i^{(0)}-2t_im\binom{m-1}{k-2},\\
     \mu_i\,t_i&=-2\binom{m-1}{k-2}+t_i(-1-2r_1(k-1)).
 \end{align*}
 On solving the system, we get
$$t_i=\frac{1+2r_1(k-1)+\mu_i^{(0)}\pm \sqrt{\bigl(1+2r_1(k-1)+\mu_i^{(0)}\bigr)^2+16m\mbinom{m-1}{k-2}^2}}{4m\mbinom{m-1}{k-2}}=\frac{-2\mbinom{m-1}{k-2}}{1+2r_1(k-1)+\mu_i}$$
and 
$$\mu_i=\frac{\mu_i^{(0)}-1-2r_1(k-1)\pm \sqrt{\bigl(1+2r_1(k-1)+\mu_i^{(0)}\bigr)^2+16m\mbinom{m-1}{k-2}^2}}{2},\mu_i^{(0)}\neq n-1-2r_0(k-1).$$
 Now we can see that, for $i\in\mathbb{Z}^{[1,n]} $
 $$S(G^*)\begin{bmatrix}
     \mathbf{0}_{1,n}\\
     \mathbf{y}_j\otimes e_i
 \end{bmatrix}=
 \mu_j^{(1)}\begin{bmatrix}
     \mathbf{0}_{1,n}\\
     \mathbf{y}_j\otimes e_i
 \end{bmatrix},$$
 where $\mathbf{y}_j$'s are the eigenvectors corresponding to the eigenvalues $\mu_j^{(1)}$'s $\,(\mu_j^{(1)}\neq m-1-2r_1(k-1)$ and $j\in\mathbb{Z}^{[1,m]}$). Thus $\mu_j^{(1)}$'s are eigenvalues of $S(G^*)$ with multiplicity $n$.\\
 
 We obtain the remaining eigenvalues of $S(G^*)$ from the eigenvalues of the quotient matrix $Q$ of $S(G^*)$,
 $$Q=\begin{bmatrix}
     n-1-2r_0(k-1) & m\bigl(n-2\binom{m-1}{k-2}\bigr)\\
     n-2\binom{m-1}{k-2} & mn-1-2r_1(k-1)
 \end{bmatrix}.$$
Thus the remaining two eigenvalues of  $S(G^*)$ can be obtained from the characteristic polynomial of $Q$ $(=\mu^2+\bigl(2-n(m+1)+2(r_0+r_1)(k-1)\bigr)\mu+(n-1-2r_0(k-1))(mn-1-2r_1(k-1))-m\bigl(n-2\binom{m-1}{k-2}\bigr)=0).$ Hence the theorem.
\end{proof}
\subsection{Coronal of the hypergraph}
Now we introduce a new invariant called coronal of a matrix, which is helpful for determining the spectrum of the corona of two hypergraphs.
\begin{defn}\label{def:coronal}
    Let $G^*$ be a $k$-uniform hypergraph on $n$ vertices with adjacency matrix $A$. If $A-\lambda I$ is invertible, then the coronal $\rchi_{A(G^*)}(\lambda)$ is defined to be the sum of entries of $(A-\lambda I)^{-1}.$ This can be written as
    $$\rchi_{A(G^*)}(\lambda)=J_{1,n}(A-\lambda I)^{-1}J_{n,1}.$$
\end{defn}
By applying Definition \ref{def:coronal}, we determine the characteristic polynomial of the corona product of two k-uniform hypergraphs in a more computationally easier manner.
\begin{thm}
Let $G_0^*$ and $G_1^*$ be $k$-uniform hypergraphs of order $n$ and $m$, respectively. Then, the characteristic polynomial of $G^* = G_0^* \odot_1^n G_1^*$ is 
    $$ P_{A(G^*)}(\lambda)=\Bigl(P_{A(G_1^*)}(\lambda)\Bigr)^nP_{A(G_0^*)}\Bigl(\lambda+\binom{m-1}{k-2}^2 \rchi_{A(G_1^*)}(\lambda) \Bigr).$$
\end{thm}
\begin{proof}
The characteristic polynomial of $G^*$ is,
     $$P_{A(G^*)}(\lambda)=\det\begin{pmatrix}
          A(G_0^*)-\lambda I_n & \binom{m-1}{k-2}J_{1,m}\otimes I_n\\
        \binom{m-1}{k-2}J_{m,1}\otimes I_n & (A(G_1^*)-\lambda I_m) \otimes I_n
     \end{pmatrix}.
     $$
     From Lemma \ref{detblock}, we get
     \begin{align*}
          P_{A(G^*)}(\lambda)&=\det\Bigl(A(G_0^*)-\lambda I_n-\mbinom{m-1}{k-2}^2 (J_{1,m}\otimes I_n) ((A(G_1^*)-\lambda I_m) \otimes I_n)^{-1}(J_{m,1}\otimes I_n) \Bigr)\\
          &\,\hspace{10cm}\det\Bigl((A(G_1^*)-\lambda I_m) \otimes I_n\Bigr)\\
           &= \det\Bigl(A(G_0^*)-\Bigl(\mbinom{m-1}{k-2}^2 \rchi_{A(G_1^*)}(\lambda) +\lambda \Bigr) I_n\Bigr)\det\Bigl(A(G_1^*)-\lambda I_m\Bigr)^n\\
           &= \Bigl(P_{A(G^*_0)}\Bigl(\mbinom{m-1}{k-2}^2 \rchi_{A(G_1^*)}(\lambda) +\lambda \Bigr)\Bigr)\bigl(P_{A(G^*_1)}(\lambda)\bigr)^n.
     \end{align*}
    
\end{proof}
\begin{corlly}
    Let $G_0^*$ and $H_0^*$ are non-isomorphic cospectral $k$-uniform hypergraph, and $G_1^*$ be any $k$-uniform hypergraph. Then $G_0^* \odot_1^n G_1^* $ and $H_0^* \odot_1^n G_1^* $ are non-isomorphic cospectral $k$-uniform hypergraphs.
\end{corlly}
By finding the corona product with $G_1^*$ repeatedly, we get infinite pair of adjacency-cospectral hypergraphs. Clearly $G_1^* \odot_1^n G_0^* $ and $G_1^* \odot_1^n H_0^* $ need not be cospectral $k$-uniform hypergraph. In addition, if $\rchi_{A(G_0^*)}(\lambda) =\rchi_{A(H_0^*)}(\lambda)$, then $G_1^* \odot_1^n G_0^* $ and $G_1^* \odot_1^n H_0^* $ are non-isomorphic cospectral hypergraphs. 

Next, we determine the coronal of  $(k,r)$-regular hypergraph, which can be  used to find the spectrum of the corona of two hypergraphs.

\begin{figure}[H]
\hspace{.5cm}
\begin{minipage}{0.2\textwidth}
\centering
\begin{tikzpicture}[scale=0.3]
\draw[ draw= black,rotate=270] (0,12) ellipse (1cm and 3cm);
\filldraw[fill=black](12,0)circle(0.1cm);

\filldraw[fill=black](10.5,0)circle(0.1cm);

\filldraw[fill=black](13.5,0)circle(0.1cm);
\end{tikzpicture}
                \vspace{2.15cm}
\subcaption{$K_3^3$}
\end{minipage}
\hspace{.5cm}
\begin{minipage}{0.2\textwidth}
\centering
\begin{tikzpicture}[scale=0.3]

\begin{scope} [fill opacity = 0]
\draw[fill=white, draw = black] (-1.5,0) circle (3);
    \draw[fill=white, draw = black] (1.5,0) circle (3);
\end{scope}
\filldraw[fill=black](0,1)circle(0.1cm);

 \filldraw[fill=black](0,-1)circle(0.1cm);

\filldraw[fill=black](-2.5,0)circle(0.1cm);

 \filldraw[fill=black](2.5,0)circle(0.1cm);

\end{tikzpicture}
                \vspace{1cm}
\subcaption{$G_0^*$}
\end{minipage}
\hspace{.5cm}
\begin{minipage}{0.5\textwidth}
\centering
\begin{tikzpicture}[scale=0.3]

\draw (36,-8) ellipse (3cm and 1cm);
\filldraw[fill=black](46,0)circle(0.1cm);
\filldraw[fill=black](46,1)circle(0.1cm);
\filldraw[fill=black](46,-1)circle(0.1cm);
\filldraw[fill=black](26,0)circle(0.1cm);
\filldraw[fill=black](26,1)circle(0.1cm);
\filldraw[fill=black](26,-1)circle(0.1cm);

\draw[dashed](26,3)--(32.5,1);
\draw[dashed](26,-3)--(32.5,-1);
\draw[dashed] (32.5,-1) to [out=20,in=-20, looseness=2.5] (32.5,1);

\draw[dashed](46,3)--(39.5,1);
\draw[dashed](46,-3)--(39.5,-1);
\draw[dashed] (39.5,1) to [out=200,in=160, looseness=2.5] (39.5,-1);
\draw[dashed](39,8)--(37,1.7);
\draw[dashed](33,8)--(35,1.7);
\draw[dashed] (35,1.7) to [out=290,in=250, looseness=2.5] (37,1.7);

\draw[dashed](33,-8)--(35,-1.7);
\draw[dashed](39,-8)--(37,-1.7);
\draw[dashed] (35,-1.7) to [out=70,in=110, looseness=2.5] (37,-1.7);

\draw (46,0) ellipse (1cm and 3cm);
\draw (26,0) ellipse (1cm and 3cm);
\filldraw[fill=black](36,-8)circle(0.1cm);

\filldraw[fill=black](34.5,-8)circle(0.1cm);

\filldraw[fill=black](37.5,-8)circle(0.1cm);

\draw (36,8) ellipse (3cm and 1cm);
\filldraw[fill=black](36,8)circle(0.1cm);

\filldraw[fill=black](34.5,8)circle(0.1cm);

\filldraw[fill=black](37.5,8)circle(0.1cm);

\draw(37.5,0)circle(3cm);
\draw(34.5,0)circle(3cm);
\filldraw[fill=black](36,1)circle(0.1cm);

 \filldraw[fill=black](36,-1)circle(0.1cm);

\filldraw[fill=black](39,0)circle(0.1cm);

 \filldraw[fill=black](33,0)circle(0.1cm);
 \node at (22,0){$G^*:$};
\end{tikzpicture}
\centering
\subcaption{ $G^*=G_0^* \odot_1^n K_3^3$}
\end{minipage}
  \caption{Example of Corona product of two hypergraphs}
  \label{fig:1}
\end{figure}
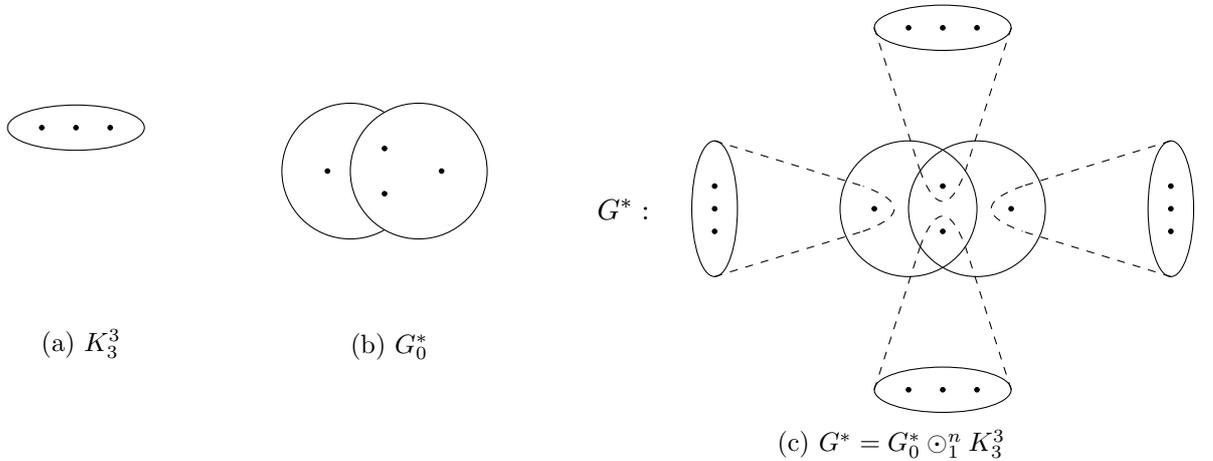

\begin{thm}
    Let $G_1^*$ be a $(k,r_1)$-regular hypergraph on $m$ vertices. Then coronal of $A(G_1^*)$ is given by,
    $$\rchi_{A(G^*_1)}(\lambda)=\frac{m}{r_1(k-1)-\lambda}.$$
    \begin{proof}
        As a consequence of regularity of $A(G_1^*)$, it follows that $(A(G_1^*)-\lambda I_m)J_{m,1}=\bigl(r(m-1)-\lambda\bigr)J_{m,1}$ and also $J_{m,1}=\displaystyle\frac{(A(G_1^*)-\lambda I_m)J_{m,1}}{\bigl(r(m-1)-\lambda\bigr)}$. Hence the coronal of $A(G_1^*)$ is,
         $$\rchi_{A(G^*_1)}(\lambda)=J_{1,m}(A(G_1^*)-\lambda I_m)^{-1}J_{m,1}=\frac{J_{1,m}J_{m,1}}{r_1(k-1)-\lambda}=\frac{m}{r_1(k-1)-\lambda}.$$
    \end{proof}
    
    \end{thm}
    Using the coronal of $(k,r)$-regular hypergraph $G_1^*$, we get the characteristic polynomial of $G^*=G_0^*\odot_1^n G_1^*$ as follows
    \begin{equation}\label{coronal_(k,r)}
           P_{A(G^*)}(\lambda)=\Bigl(P_{A(G_1^*)}(\lambda)\Bigr)^nP_{A(G_0^*)}\Bigl(\lambda+ \frac{m}{r(k-1)-\lambda} \binom{m-1}{k-2}^2\Bigr). 
    \end{equation}
    We can easily find the spectrum of $G^*$ from (\ref{coronal_(k,r)}), if the spectrum of $G_0^*$ and $G_1^*$ are known. The following example gives the spectrum of $G^*$ if $G_1^*=K_3^3$ and $G_0^*$ as in Figure \ref{fig:1}(b).

\begin{exmp}\label{coronalexmp}
    \textup{ Let $G_0^*$ and $K_3^3$ be 3-uniform hypergraphs as given in Figure \ref{fig:1} and  $G^*=G_0^* \odot_1^n K_3^3$. The characteristic polynomial of $A(G_0^*)$ is  $P_{A(G^*_0)}(\lambda)=\lambda^4 - 8\lambda^2 - 8\lambda$, $P_{A(K_3^3)}(\lambda)=- \lambda^3 + 3\lambda + 2$ and the coronal of $K_3^3$ is  $\rchi_{K_3^3}(\lambda)=\frac{-3}{\lambda - 2}$. Here $n=4$, $m=3$ and $k=3$, we have }
  \begin{align*}
      P_{A(G^*)}(\lambda)&=(- \lambda^3 + 3\lambda + 2)^4\biggl(\Bigl(\lambda-\frac{12}{\lambda - 2}\Bigr)^4 - 8\Bigl(\lambda-\frac{12}{\lambda - 2}\Bigr)^2 - 8\Bigl(\lambda-\frac{12}{\lambda - 2}\Bigr)\biggr)\\
       &=(\lambda+1)^8(\lambda^2-16)(\lambda^2-2\lambda-12)(\lambda^4-6\lambda^3-16\lambda^2+80\lambda+80).
  \end{align*}
  \textup{Then the spectrum of $G^*$ is,}
  $$\sigma_A({G^*})=
    \begin{pmatrix}
    -1&16&-16&\displaystyle\frac{829}{180}&\displaystyle \frac{-3567}{1369} & \displaystyle\frac{10046}{1637} &\displaystyle\frac{1106}{263}  &\displaystyle\frac{-2877}{836}& \displaystyle\frac{-1525}{1693} \\
              8 &1&1&   1 & 1 & 1 & 1 & 1 & 1
    \end{pmatrix}$$
\end{exmp}

\section{Spectrum of the corona hypergraph $G_0^{*(m)}$}
In this section, we define the corona hypergraph based on the definition of the corona product of two hypergraphs. In addition,  we determine the spectrum of the corona hypergraph.
\begin{definition}
    Let $G_0^*$ be a $k$-uniform hypergraph on $n$ vertices. Then the corona hypergraph $G_0^{*(m)}$ with respect to the base hypergraph $G_0^*=G_0^{*(1)}$ is obtained by,
    $$G_0 ^{*(m)}=G_0^{*(m-1)} \odot_1^n G_0^{*(1)}. $$
\end{definition}
 \begin{figure}[H]
 \centering
\begin{tikzpicture}[scale=.3]

\draw[rotate=45] (19,-11) ellipse (1.6cm and 0.7cm);
\filldraw [rotate=45] [fill=black](18,-11)circle(0.1cm); 
\filldraw [rotate=45] [fill=black](19,-11)circle(0.1cm); 
\filldraw [rotate=45] [fill=black](20,-11)circle(0.1cm); 

\draw[dashed](22.35,6.7)--(30.5,0.6);
\draw[dashed](20.1,4.5)--(30,-0.4);
\draw [dashed] (30,-0.4) .. controls (31,-0.5) and (31,-0.5) .. (30.5,0.6);

\draw[rotate=145] (-24.5,-8.5) ellipse (1.6cm and 0.7cm);
\filldraw [rotate=145] [fill=black](-23.5,-8.5)circle(0.1cm); 
\filldraw [rotate=145] [fill=black](-24.5,-8.5)circle(0.1cm); 
\filldraw [rotate=145] [fill=black](-25.5,-8.5)circle(0.1cm);

\draw[dashed](23.75,-6)--(30,0.5);
\draw[dashed](26.22,-8.1)--(31,-0.6);
\draw [dashed] (30,0.5) .. controls (31,0.5) and (31,-0.5)..(31,-0.6);

\draw[rotate=135] (-32,-40) ellipse (1.6cm and 0.7cm);
\filldraw [rotate=135] [fill=black](-31,-40)circle(0.1cm); 
\filldraw [rotate=135] [fill=black](-32,-40)circle(0.1cm); 
\filldraw [rotate=135] [fill=black](-33,-40)circle(0.1cm); 

\draw[dashed](52,4.5)--(42.5,-0.4);
\draw[dashed](49.75,6.7)--(42,0.9);
\draw [dashed] (42,0.9) .. controls (41,0.5) and (41,-0.5) .. (42.5,-0.4);

\draw (36,0) ellipse (6.5cm and 1.5cm);
\filldraw[fill=black](41.5,0)circle(0.1cm); 
\filldraw[fill=black](36,0)circle(0.1cm); 
\filldraw[fill=black](30.5,0)circle(0.1cm);

 \draw (30.5,4) ellipse (3cm and 1cm);
\filldraw[fill=black](28,4)circle(0.1cm); 
\filldraw[fill=black](30.5,4)circle(0.1cm);
\filldraw[fill=black](33,4)circle(0.1cm);
\draw[dashed](33.5,4)--(31.5,1);
\draw[dashed](27.5,4)--(29.5,1);
\draw [dashed] (29.5,1) .. controls (30.5,-0.8) and (30.5,-0.8) .. (31.5,1);
\draw (30.5,8) ellipse (1.6cm and 0.7cm);
\filldraw[fill=black](29.4,8)circle(0.1cm); 
\filldraw[fill=black](30.5,8)circle(0.1cm); 
\filldraw[fill=black](31.6,8)circle(0.1cm);

 \draw[dashed](28.9,8)--(30,4.5);
\draw[dashed](32.1,8)--(31,4.5);
\draw [dashed] (30,4.5) .. controls (30.5,3.3) and (30.5,3.3) .. (31,4.5);

\draw (27.1,8) ellipse (1.6cm and 0.7cm);
\filldraw[fill=black](26,8)circle(0.1cm); 
\filldraw[fill=black](27.1,8)circle(0.1cm); 
\filldraw[fill=black](28.2,8)circle(0.1cm);

 \draw[dashed](25.5,8)--(27.4,4.5);
\draw[dashed](28.7,8)--(28.4,4.5);
\draw [dashed] (27.4,4.5) .. controls (28,3.4) and (28.3,3.4) .. (28.4,4.5);

\draw (33.9,8) ellipse (1.6cm and 0.7cm);
\filldraw[fill=black](32.8,8)circle(0.1cm); 
\filldraw[fill=black](33.9,8)circle(0.1cm); 
\filldraw[fill=black](35,8)circle(0.1cm);

\draw[dashed](32.3,8)--(32.5,4.5);
\draw[dashed](35.5,8)--(33.5,4.5);
\draw [dashed] (32.5,4.5) .. controls (32.7,3.4) and (33,3.4) .. (33.5,4.5);

 \draw (41.5,4) ellipse (3cm and 1cm);
\filldraw[fill=black](39,4)circle(0.1cm); 
\filldraw[fill=black](41.5,4)circle(0.1cm);
\filldraw[fill=black](44,4)circle(0.1cm);

\draw[dashed](44.5,4)--(42.5,1);
\draw[dashed](38.5,4)--(40.5,1);
\draw [dashed] (40.5,1) .. controls (41.5,-0.8) and (41.5,-0.8) .. (42.5,1);

\draw (41.5,8) ellipse (1.6cm and 0.7cm);
\filldraw[fill=black](40.4,8)circle(0.1cm); 
\filldraw[fill=black](41.5,8)circle(0.1cm); 
\filldraw[fill=black](42.6,8)circle(0.1cm);

\draw[dashed](39.9,8)--(41,4.5);
\draw[dashed](43.1,8)--(42,4.5);
\draw [dashed] (41,4.5) .. controls (41.5,3.2) and (41.5,3.2) .. (42,4.5);

\draw (38.1,8) ellipse (1.6cm and 0.7cm);
\filldraw[fill=black](37,8)circle(0.1cm); 
\filldraw[fill=black](38.1,8)circle(0.1cm); 
\filldraw[fill=black](39.2,8)circle(0.1cm);

\draw[dashed](36.5,8)--(38.4,4.4);
\draw[dashed](39.7,8)--(39.4,4.4);
\draw [dashed] (38.4,4.4) .. controls (39,3.5) and (39,3.5) .. (39.4,4.4);

\draw (44.9,8) ellipse (1.6cm and 0.7cm);
\filldraw[fill=black](43.8,8)circle(0.1cm); 
\filldraw[fill=black](44.9,8)circle(0.1cm); 
\filldraw[fill=black](46,8)circle(0.1cm);

\draw[dashed](43.3,8)--(43.5,4.5);
\draw[dashed](46.5,8)--(44.5,4.5);
\draw [dashed] (43.5,4.5) .. controls (44,3.3) and (44,3.3) .. (44.5,4.5);

\draw (36,-4) ellipse (3cm and 1cm);
\filldraw[fill=black](38,-4)circle(0.1cm); 
\filldraw[fill=black](36,-4)circle(0.1cm); 
\filldraw[fill=black](34,-4)circle(0.1cm);

\draw[dashed](33,-4)--(35,-1);
\draw[dashed](39,-4)--(37,-1);
\draw[dashed] (35,-1) to [out=70,in=110, looseness=2.5] (37,-1);

\draw (36,-8) ellipse (1.6cm and 0.7cm);
\filldraw[fill=black](35.3,-8)circle(0.1cm); 
\filldraw[fill=black](36,-8)circle(0.1cm); 
\filldraw[fill=black](36.7,-8)circle(0.1cm);
 
\draw[dashed](34.4,-8)--(35.3,-5);
\draw[dashed](37.6,-8)--(36.7,-5);
\draw[dashed] (35.3,-4.6) to [out=70,in=110, looseness=2.5] (36.7,-4.6);

\draw (32,-8) ellipse (1.6cm and 0.7cm);
\filldraw[fill=black](31.3,-8)circle(0.1cm); 
\filldraw[fill=black](32,-8)circle(0.1cm); 
\filldraw[fill=black](32.7,-8)circle(0.1cm);

 \draw[dashed](30.4,-8)--(33,-4.5);
\draw[dashed](33.6,-8)--(34.5,-4.5);
\draw [dashed] (33,-4.5) .. controls (33.6,-3.2) and (34.5,-3.2) .. (34.5,-4.5);

\draw (40,-8) ellipse (1.6cm and 0.7cm);
\filldraw[fill=black](39.3,-8)circle(0.1cm); 
\filldraw[fill=black](40,-8)circle(0.1cm); 
\filldraw[fill=black](40.7,-8)circle(0.1cm);
 
\draw[dashed](38.4,-8)--(37.5,-4.5);
\draw[dashed](41.6,-8)--(39,-4.5);
\draw [dashed] (37.5,-4.5) .. controls (37.6,-3.2) and (38.7,-3.2) .. (39,-4.5);

 \end{tikzpicture}
 \caption{$K_3^{3(3)}$}
\label{fig:2}
\end{figure}
Figure \ref{fig:2} illustrates the corona hypergraph $G_0^{*(m)}$, when $G_0^{*}=K_m^m$.\\
Following are some results on the corona hypergraphs.
\begin{enumerate}
    \item The order of $G_0^{*(m)}$ is
            $n(n+1)^{m-1}.$
    \item The size of  $G_0^{*(m)}$ is  $\tilde{e}(n+1)^{(m-1)}+\binom{n}{k-1}\bigl((n+1)^{m-2}-1\bigr),$ where $\tilde{e}$ is the size of $G_0^{*}$.
    \item If $G_0^{*}$ is connected, then the corona hypergraph $G_0^{*(m)}$ is also connected.
\end{enumerate}

Now we obtain the adjacency spectrum of the corona hypergraph $G_0^{*(m)}$ in terms of the eigenvalues of $G_0^*$ by applying Theorem \ref{corona2hygrph}.\\ 

For any $(k,r)$-regular hypergraph of order $n$, we denote
\begin{equation*}
    \phi(\lambda)=\frac{\lambda+r(k-1)\pm \sqrt{\bigl( \lambda-r(k-1))^2+4n\binom{n-1}{k-2}^2\bigr)}}{2},
\end{equation*}
 $\phi^0(\lambda)=\lambda$ and $\phi^m(\lambda)=\phi\bigl(\phi^{m-1}(\lambda)\bigr)$.
 \begin{theorem}
     Let $G_0^*$ be a $(k,r)$ regular hypergraph and its spectrum be
     $$\sigma_A(G_0^*)=\begin{pmatrix}
         \lambda_1=r(k-1) & \lambda_2 & \lambda_3 &\cdots&\lambda_n\\
         1 & 1 &1&\cdots&1
     \end{pmatrix}.$$
     Then the spectrum of corona hypergraph $G_0^{*(m)}$ is,
     \begin{equation*}
         \sigma_A(G_0^{*(m)})=\displaystyle\begin{pmatrix}
         \phi^{m-1}(\lambda_i)_{i\in\mathbb{Z}^{[1,n]}} & \phi^{m-2}(\lambda_i)_{i\in\mathbb{Z}^{[2,n]}}&\cdots& \phi^j(\lambda_i)_{i\in\mathbb{Z}^{[2,n]}} &\cdots&\phi^0(\lambda_i)_{i\in\mathbb{Z}^{[2,n]}}\\
         1 & n &\cdots &n(n+1)^{m-j-2}&\cdots&n(n+1)^{m-2}
     \end{pmatrix}.
     \end{equation*}
 \end{theorem}
 \begin{proof}
    To prove this we use mathematical induction on $m$. From Theorem \ref{corona2hygrph}, for $i\in\mathbb{Z}^{[1,n]}, j\in\mathbb{Z}^{[2,n]}$ the spectrum of  $G_0^{*(2)}=G_0^*\odot G_0^*$ is given by
    $$\sigma_A({G_0^{*(2)}})=
    \begin{pmatrix}
     \frac{\lambda_i+r(k-1)\pm \sqrt{\left(\lambda_i-r(k-1)\right)^2+4n\mbinom{n-1}{k-2}^2}}{2}  &\lambda_j \\
                 1 &n
    \end{pmatrix}=\displaystyle\begin{pmatrix}
         \phi^1(\lambda_i)_{i\in\mathbb{Z}^{[1,n]}} & \phi^0(\lambda_i)_{i\in\mathbb{Z}^{[2,n]}}\\
         1 & n 
     \end{pmatrix}.$$
    Since $\phi^1(\lambda_i)_{i\in\mathbb{Z}^{[1,n]}}$ gives $2n$ eigenvalues and  $\phi^0(\lambda_i)_{i\in\mathbb{Z}^{[2,n]}}$ gives the rest of the $n(n-1)$ eigenvalues of $G_0^{*(2)}$, we get all the $n(n+1)$ eigenvalues of $G_0^{*(2)}$. \\

    Let the result holds for $G_0^{*(m)}(m>2).$ Next we prove for $G_0^{*(m+1)}=G_0^{*(m)}\odot_1^n G_0^*$. By Theorem \ref{corona2hygrph}, we have
    \begin{align}\label{G*m.1}
         \sigma_A(G_0^{*(m+1)})&=\begin{pmatrix}
        \phi(\sigma_A(G_0^{*(m)})) &\phi^0(\lambda_i)_{i\in\mathbb{Z}^{[2,n]}}\\
        1 & n
    \end{pmatrix}\notag \\
    &=\begin{pmatrix}
        \phi\bigl(\phi^{(m-1)}(\lambda_i)_{i\in\mathbb{Z}^{[1,n]}}\bigr) & \phi\bigl(\phi^{(m-2)}(\lambda_i)_{i\in\mathbb{Z}^{[2,n]}}\bigr) &\cdots &\phi\bigl(\phi^0(\lambda_i)_{i\in\mathbb{Z}^{[2,n]}}\bigr) & \phi^0(\lambda_i)_{i\in\mathbb{Z}^{[2,n]}}\\
        1 & n & \cdots &n(n+1)^{m-2} &n(n+1)^{m-1}
    \end{pmatrix}.
    \end{align}
  It is clear that $\phi^{m}(\lambda_i)_{i\in\mathbb{Z}^{[1,n]}}$ gives $2^m n$ eigenvalues and  $\phi^j(\lambda_i)_{i\in\mathbb{Z}^{[2,n]}}$ gives  $2^j n(n-1)(n+1)^{m-j-1}$ eigenvalues of $G_0^{*(m+1)}$. Therefore, from (\ref{G*m.1}) and  Lemma \ref{noofeig}  we obtain all the $n(n+1)^m$ eigenvalues of $G_0^{*(m+1)}$.  Hence the theorem.
 \end{proof}

 \section{Construction of Seidel-cospectral hypergraphs}
 
Sarkar and Banerjee \cite{Sarkar2020} introduced a method for constructing cospectral $k$-uniform hypergraphs based on the vertex corona of hypergraphs. Later, in \cite{abiad2022}, various techniques for generating non-isomorphic adjacency-cospectral hypergraphs were discussed. 
    
In this section, we focus on constructing non-isomorphic Seidel-cospectral $k$-uniform hypergraphs by extending the method of signed M-GM switching \cite{Belardo2021} to find a family of switching equivalent hypergraphs. \\

Let $G^*$ be a $k$-uniform hypergraph with vertex set $V$ and hyperedge set $E$, and $ \{U_1, U_2, U_3,\ldots, U_{2t}, U\},$ $t\in \mathbb{N}$ be the partition of $V$. Then, proper labelling of the vertices, enables us to express the Seidel matrix of $G^*$ as follows:
   \begin{equation*}
       S(G^*)=\begin{bmatrix}
       S_{1,1} & S_{1,2} & S_{1,3}&\cdots&S_{1,2t}& S_1\\
       S_{2,1} & S_{2,2} & S_{2,3}&\cdots&S_{2,2t}& S_2\\
       S_{3,1} & S_{3,2} & S_{3,3}&\cdots&S_{3,2t}& S_3\\
       \vdots & \vdots & \vdots&\ddots&\vdots& \vdots\\
       S_{2t,1} & S_{2t,2} & S_{2t,3}&\cdots&S_{2t,2t}& S_{2t}\\
       S_1^T & S_2^T & S_3^T&\cdots & S_{2t}^T & S_0
            \end{bmatrix}.
    \end{equation*}
    \begin{theorem}
        Let $G^*=(V, E)$ be a $k$- uniform hypergraph and $ \{U_1, U_2, U_3,\ldots, U_{2t}, U\}$ be the partition of V  that  satisfies the following conditions:
        \begin{itemize}
            \item $|U_i|=m,\;m\in \mathbb{Z}^+, ~0< i \leq 2t$ and $|U|=k-1$. .
            \item For any odd $i$, $N(U)\cap(U_i \cup U_{i+1})= U_i$, where $N(U)=\{v \in V: v\cup U\in E\} $.
            \item For any odd $i,j~(1\leq i,j \leq 2t)$, 
            \begin{itemize}
                \item[] for $1\leq s \leq m,~ \sum_{r=1}^m \left((S_{i,j})_{rs}-(S_{i,j+1})_{rs}\right)= \sum_{r=1}^m((S_{i+1,j+1})_{rs}-(S_{i+1,j})_{rs})=l$, 
                \item[] for $1\leq r \leq m,~ \sum_{s=1}^m \left((S_{i,j})_{rs}-(S_{i+1,j})_{rs}\right)= \sum_{s=1}^m((S_{i+1,j+1})_{rs}-(S_{i,j+1})_{rs})=l$.
            \end{itemize}   
        \end{itemize}
Then a non-isomorphic Seidel cospectral mate $k$-uniform hypergraph $H^*$ of $G^*$ forms by replacing the edge $\{u_i\}\cup U$,  $\forall u_i\in U_i$ with the new edges $\{u_{i+1}\} \cup U$, $u_{i+1}\in U_{i+1}$.

    \end{theorem}
    \begin{proof}
        Consider an orthogonal block matrix $P$  
        \begin{equation*}
       P=\begin{bmatrix}
       P_1 & P_2 & \mathbf{O}&\mathbf{O}&\cdots& \mathbf{O}\\
       P_2  & P_1 & \mathbf{O}&\mathbf{O}&\cdots&\mathbf{O}\\
       \mathbf{O} &\mathbf{O} &P_1&P_2 &\cdots&\mathbf{O}\\
       \mathbf{O}& \mathbf{O} &P_2 &P_1&\cdots& \mathbf{O}\\
       \vdots & \vdots & \vdots& \vdots& \ddots& \vdots\\
       \mathbf{O}& \mathbf{O} & \mathbf{O}&\mathbf{O}&\cdots & I_{k-1}
            \end{bmatrix}\;\text{and}\; 
        \tilde{S}=\begin{bmatrix}
       \tilde{S}_{1,1} & \tilde{S}_{1,2} & \cdots &\tilde{S}_{1,2t}&\tilde{S}_1\\
       \tilde{S}_{2,1} & \tilde{S}_{2,2} & \cdots &\tilde{S}_{2,2t}&\tilde{S}_2\\
       \vdots & \vdots & \ddots & \vdots & \vdots\\
       \tilde{S}_{2t,1} & \tilde{S}_{2t,2} & \cdots & \tilde{S}_{2t,2t}&\tilde{S}_2t\\
       \tilde{S}_{1} & \tilde{S}_{2} & \cdots & \tilde{S}_{2t}&\tilde{S}_0
            \end{bmatrix}
    \end{equation*}
 such that $PSP=\tilde{S}$, where $P_1=I_m-\frac{1}{m}J_m$ and $P_2=\frac{1}{m}J_m$. Therefore, for any odd $i,~j$
 \begin{equation*}
   \begin{bmatrix}
       P_1 & P_2 \\
       P_2  & P_1  
    \end{bmatrix}
    \begin{bmatrix}
        S_{i,j} & S_{i,j+1} \\
       S_{i+1,j} & S_{i+1,j+1} 
    \end{bmatrix}
     \begin{bmatrix}
       P_1 & P_2 \\
       P_2  & P_1  
    \end{bmatrix}=
    \begin{bmatrix}
       \tilde{ S}_{i,j} & \tilde{S}_{i,j+1} \\
      \tilde{S}_{i+1,j} & \tilde{S}_{i+1,j+1} 
    \end{bmatrix}.
 \end{equation*}
Now we prove that $\tilde{ S}_{i,j} =S_{i,j}$
    \begin{align*}
        \tilde{ S}_{i,j}& = \bigl(P_1S_{i,j}+P_2 S_{i+1,j}\bigr)P_1+ \bigl(P_1 S_{i,j+1}+P_2 S_{i+1,j+1}\bigr)P_2      \\ 
                        \begin{split}
                            = S_{i,j}+ \frac{1}{m}\Bigl(-J_m \bigl(S_{i,j}-S_{i+1,j}\bigr)- \bigl(S_{i,j}-S_{i,j+1}\bigr) J_m\Bigr) +\frac{1}{m^2}\Bigl(J_m  \bigl(S_{i,j}-S_{i,j+1}\bigr)J_m\\+J_m \bigl(S_{i+1,j+1}
                        -S_{i+1,j}\bigr)J_m\Bigr) 
                        \end{split}\\
            &=S_{i,j}+ \frac{1}{m}\Bigl(-l\, J_m - l\, J_m\Bigr) +\frac{1}{m^2}\Bigl( l\,m\, J_m+l\, m\, J_m\Bigr)\\
            &=S_{i,j}.
    \end{align*}
   
    From a similar argument we obtain,
     \begin{equation*}
    \begin{bmatrix}
       \tilde{ S}_{i,j} & \tilde{S}_{i,j+1} \\
      \tilde{S}_{i+1,j} & \tilde{S}_{i+1,j+1} 
    \end{bmatrix}=
     \begin{bmatrix}
        S_{i,j} & S_{i,j+1} \\
       S_{i+1,j} & S_{i+1,j+1} 
    \end{bmatrix}.
 \end{equation*}
 Next, we will show that for an odd $i\,(0\leq i \leq 2t)$, $\tilde{ S}_{i}= S_{i+1}$ and $\tilde{ S}_{i+1}= S_{i}$ according to which all the values of Seidel matrix corresponding to $k-1$ entries of $U$ and $U_i$ is switched with $U_{i+1}$. Since $N(U)\cap(U_i \cup U_{i+1})= U_i$, we can observe that $S_i$ and $S_{i+1}$ are block matrix with all entries equal to $-1$ and  $1$ respectively. Then, $J_mS_i=-m\,J_{m,k-1}$ and$J_mS_{i+1}=m\,J_{m,k-1}$. Therefore,
 \begin{align*}
  \tilde{ S}_{i}= P_1S_i+P_2S_{i+1}=S_i+2J_{m,k-1}=S_{i+1}
\end{align*}
and
\begin{align*}
  \tilde{ S}_{i+1}= P_2S_i+P_1S_{i+1}=S_{i+1}-2J_{m,k-1}=S_{i}
\end{align*}
Hence, we can conclude that $S$ and $\tilde{S}$ are similar matrices corresponding to Seidel-cospectral hypergraphs $G^*$ and $H^*$, respectively.
    \end{proof}
\begin{example}
    Let $G^*$ b a hypergraph with vertex set $V(G^*)=\{v_0, v_1, u_1, u_2, u_3,u_4,u_5,u_6\}$ and edge set $E=\{\{u_1,u_2,u_3\},\{u_1,u_4,u_5\},\{u_2,u_5,u_6\},\{u_3,u_4,u_6\},\{v_0,v_1,u_1\},\{v_0,v_1,u_2\},\{v_0,v_1,u_3\}\}$.Consider the partition $\{U_1=\{u_1,u_2,u_3\},U_2=\{u_4,u_5,u_6\},U=\{v_0,v_1,\}\}$ of $V$. Then  $N(U)=\{u_1,u_2,u_3\}=U_1$. Now, we get a non-regular  non-isomorphic Seidel-cospectral hypergraph $H^*$ with vertex set $V$ and edge set $E/\{\{v_0,v_1,u_1\},$ $\{v_0,v_1,u_2\},\{v_0,v_1,u_3\}\}\cup\{\{v_0,v_1,u_4\},\{v_0,v_1,u_5\},\{v_0,v_1,u_6\}\}$. 
\end{example}
Next, we prove the existence of infinitely many non-isomorphic Seidel-cospectral hypergraphs for any two hypergraphs that are not isomorphic.
\begin{theorem}
    Let $G_1^*$ and $H_1^*$ be non-isomorphic cospectral $(k,r)$-regular hypergraphs of order $m$. Then $G_0^* \odot_1^n G_1^*$ and $G_0^* \odot_1^n H_1^*$ are non-regular Seidel-cospectral hypergraphs.
\end{theorem}
\begin{proof}
    The proof follows from Theorem \ref{S-cosepec_1} and Theorem \ref{S-cospec2}.
\end{proof}

\section{Conclusion}   
 In this paper, the characteristic polynomial of generalised adjacency and generalised Seidel matrix of the corona of hypergraphs are determined. Also, proposed a method to generate infinitely many non-regular and non-isomorphic Seidel-cospectral hypergraphs. In addition, estimated the adjacency spectrum, Seidel spectrum and adjacency eigenvectors of the corona of two hypergraphs. Corona hypergraphs are constructed by iteratively applying the corona product of the base hypergraph. When the base hypergraph is regular, we determined the adjacency spectrum of corona hypergraphs. Moreover, a new invariant the coronal of the adjacency matrix of a hypergraph, has been introduced to find the spectrum of the corona of hypergraphs in a more computably easier way.

\section{Declarations}
 On behalf of all authors, the corresponding author states that there is no conflict of interest.\\

\bibliographystyle{plain}
\bibliography{ref}
\end{document}